\documentclass[reqno]{amsart}
\usepackage{amsmath}
\usepackage{amssymb}
\usepackage{amsthm}
\usepackage{cases}
\usepackage{mathtools}
\mathtoolsset{showonlyrefs=true}
\usepackage{mathrsfs}
\usepackage{hyperref}

\title{A vertex operator reformulation of the Kanade-Russell conjecture modulo 9}

\author{Shunsuke Tsuchioka}
\address{Department of Mathematical and Computing Sciences, Tokyo Institute of Technology, Tokyo 152-8551, Japan}
\email{tshun@kurims.kyoto-u.ac.jp}

\date{Jun 12, 2023}

\usepackage{tikz}

\def\node#1#2{\overset{#1}{\underset{#2}{\circ}}}

\def\ver#1#2{\overset{{\llap{$\scriptstyle#1$}\displaystyle\circ{\rlap{$\scriptstyle#2$}}}}{\scriptstyle\vert}}

\usetikzlibrary{arrows}
\tikzstyle{every picture}+=[remember picture]
\tikzstyle{na} = [baseline=-.5ex]
\tikzstyle{mine}= [arrows={angle 90}-{angle 90},thick]

\def\Llleftarrow{%
\lower2pt\hbox{\begingroup
\tikz
\draw[shorten >=0pt,shorten <=0pt] (0,3pt) -- ++(-1em,0) (0,1pt) -- ++(-1em-1pt,0) (0,-1pt) -- ++(-1em-1pt,0) (0,-3pt) -- ++(-1em,0) (-1em+1pt,5pt) to[out=-105,in=45] (-1em-2pt,0) to[out=-45,in=105] (-1em+1pt,-5pt);
\endgroup}
}

\newtheorem{Thm}{Theorem}[section]
\newtheorem{Def}[Thm]{Definition}
\newtheorem{Conj}[Thm]{Conjecture}

\newtheorem{Prop}[Thm]{Proposition}

\newtheorem{Rem}[Thm]{Remark}
\newtheorem{Cor}[Thm]{Corollary}
\newtheorem{Ex}[Thm]{Example}

\newcommand{\KE}[1]{\widetilde{e}_{#1}}
\newcommand{\KF}[1]{\widetilde{f}_{#1}}

\newcommand{\GEE}{\mathfrak{g}}
\newcommand{\GE}{\widetilde{\mathfrak{g}}}
\newcommand{\GH}{\mathfrak{t}}
\newcommand{\GA}{\mathfrak{a}}
\newcommand{\GAA}{\widetilde{\mathfrak{a}}}
\newcommand{\GAAA}{\widehat{\mathfrak{a}}}

\newcommand{\PC}{\mathcal{C}}
\newcommand{\PD}{\mathcal{D}}
\newcommand{\ZERO}{\boldsymbol{0}}
\newcommand{\TOP}{\emptyset}

\newcommand{\PT}{\stackrel{\mathsf{PT}}{\sim}}

\DeclareMathOperator{\SEQ}{\mathsf{Seq}}

\DeclareMathOperator{\PAR}{\mathsf{Par}}

\DeclareMathOperator{\IND}{\mathsf{Ind}}

\DeclareMathOperator{\ID}{id}

\DeclareMathOperator{\SUPP}{\mathsf{Supp}}

\DeclareMathOperator{\TT}{\mathsf{T^{=}}}
\DeclareMathOperator{\T}{\mathsf{T}}
\DeclareMathOperator{\Z}{\mathsf{Z}}
\DeclareMathOperator{\X}{\mathsf{X}}

\DeclareMathOperator{\END}{\mathsf{End}}

\DeclareMathOperator{\PR}{\mathsf{pr}}
\DeclareMathOperator{\PAI}{\pi^{\bullet}}
\DeclareMathOperator{\PAII}{\iota}

\newcommand{\BARR}[1]{\overline{#1}}

\newcommand{\AL}[1]{${#1}$}
\newcommand{\JO}[1]{$q_{#1}$}
\newcommand{\JOO}[1]{q_{#1}}

\newcommand{\EMPTYWORD}{\varepsilon}
\newcommand{\EMPTYPART}{\boldsymbol{\emptyset}}

\newcommand{\RR}{\mathsf{RR}}
\newcommand{\KR}{\mathsf{KR}}
\newcommand{\EL}{\mathsf{L}}
\newcommand{\LWC}{\mathsf{C}}

\newcommand{\JEI}{\mathsf{K}}
\newcommand{\AVOID}{\mathsf{avoid}}
\newcommand{\EMU}{\mathcal{M}}
\newcommand{\emu}{\mathsf{m}}

\newcommand{\EFU}{\mathsf{F}}
\newcommand{\AI}{\mathsf{I}}
\DeclareMathOperator{\SAT}{\mathsf{sat}}
\DeclareMathOperator{\A}{\mathsf{M}}
\DeclareMathOperator{\B}{\mathsf{N}}
\DeclareMathOperator{\C}{\mathsf{C}}
\DeclareMathOperator{\AD}{\mathsf{S}}
\DeclareMathOperator{\BD}{\mathsf{T}}
\DeclareMathOperator{\CD}{\mathsf{U}}
\DeclareMathOperator{\D}{\mathsf{P}}
\DeclareMathOperator{\E}{\mathsf{Q}}

\DeclareMathOperator{\QUE}{\mathsf{I}}
\DeclareMathOperator{\QUEE}{\mathsf{J}}
\DeclareMathOperator{\QUEEE}{\mathsf{K}}
\newcommand{\ICHI}[1]{\theta^{(1)}_{#1}}
\newcommand{\NI}[1]{\theta^{(2)}_{#1}}
\newcommand{\SAN}[1]{\theta^{(3)}_{#1}}
\newcommand{\SHI}[1]{\theta^{(4)}_{#1}}

\usepackage{graphicx}

\begin{document}
\maketitle

\begin{abstract}
We reformulate  the Kanade-Russell conjecture modulo 9
via the vertex operators for the level 3 standard modules of 
type $D^{(3)}_{4}$.
Along the same line, we arrive at 
three partition theorems which may be regarded as an $A^{(2)}_{4}$ analog of the conjecture.
Andrews-van Ekeren-Heluani has proven one of them, and we point out that the others are easily proven from their results.
\end{abstract}

\section{Introduction}
\subsection{The Rogers-Ramanujan identities}
The Rogers-Ramanujan (RR, for short) identities
\begin{align*}
\sum_{n\geq 0}\frac{q^{n^2}}{(q;q)_n}=\frac{1}{(q,q^4;q^5)_{\infty}},\quad
\sum_{n\geq 0}\frac{q^{n^2+n}}{(q;q)_n}=\frac{1}{(q^2,q^3;q^5)_{\infty}}  
\end{align*}
are undoubtedly one of the most famous $q$-series identities in mathematics (see ~\cite{An1,Har,Sil}).
Here, the Pochhammer symbols are defined by
\begin{align*}
  (a;q)_n = (1-a)(1-aq)\cdots (1-aq^{n-1}),\quad
  (a_1,\dots,a_k;q)_n = (a_1;q)_n\cdots (a_k;q)_n.
\end{align*}

As noted by Schur and MacMahon,
the RR identities are equivalent to the following statement, often called the RR partition theorem.
\begin{quotation}
  Let $i=1$ or $2$. For each non-negative integer $n$, the number of partitions of $n$ whose successive differences are at least two
  and whose minimum parts are at least $i$
  equals the number of partitions of $n$ whose parts are congruent to $i$ or $5-i$ modulo 5. 
\end{quotation}

Recall that a partition $\lambda=(\lambda_1,\dots,\lambda_{\ell})$ of a non-negative integer $n$ is
a weakly decreasing sequence of positive integers (called parts) whose sum $|\lambda|$ is $n$.
We denote the length $\ell$ by $\ell(\lambda)$ and put the multiplicities as $m_j(\lambda)=|\{1\leq i\leq \ell\mid \lambda_i=j\}|$.

Let $\PAR(n)$ (resp. $\PAR$) denote the set of partitions of $n$ (resp. partitions).
We say that two subsets $\PC$ and $\PD$ of $\PAR$ are partition theoretically equivalent (abbreviated to $\PC\PT\PD$)
if we have $|\PC\cap\PAR(n)|=|\PD\cap\PAR(n)|$ for $n\geq 0$.
For example, the 
RR partition theorem is briefly written as $\RR_i\PT T^{(5)}_{i,4-i}$ for $i=1,2$, where
\begin{align*}
  \RR_1 &= \{\lambda\in\PAR\mid \lambda_j-\lambda_{j+1}\geq 2 \textrm{ for $1\leq j<\ell(\lambda)$}\},\\
  \RR_2 &= \RR_1\cap\{\lambda\in\PAR\mid m_1(\lambda)=0\},\\
  T^{(N)}_{a_1,\dots,a_k} &= \{\lambda\in\PAR\mid \lambda_j\equiv a_1,\dots,a_k\!\!\pmod{N}\textrm{ for $1\leq j\leq\ell(\lambda)$}\}.
\end{align*}

\subsection{The Kanade-Russell conjecture modulo 9}
About a decade ago, Kanade-Russell found a celebrated conjectural partition theorem below (see ~\cite[\S4]{KR}), which is now well-known
in the community under the name of the Kanade-Russell conjecture (e.g., see ~\cite{CL,Kur,UZ}).

\begin{Conj}[{\cite{KR}}]\label{krconj}
  For $1\leq a\leq 3$, we have $\KR_a\PT T^{(9)}_{2^{a-1},3,6,9-2^{a-1}}$, where 
\begin{align*}
  \KR_1 &= \{\lambda\in\PAR\mid \lambda_j-\lambda_{j+2}\geq 3 \textrm{ for $1\leq j\leq \ell(\lambda)-2$ and $P(\lambda,0)$ holds}\}, \\
  \KR_2 &= \KR_1\cap\{\lambda\in\PAR\mid m_1(\lambda)=0\},\\
  \KR_3 &= \KR_1\cap\{\lambda\in\PAR\mid m_1(\lambda)=m_2(\lambda)=0\}.
\end{align*}
Moreover, the condition $P(\lambda,k)$ for a partition $\lambda$ and an integer $k$ stands for
\begin{align*}
\textrm{$\lambda_j-\lambda_{j+1}\leq 1$ implies $\lambda_j+\lambda_{j+1}\equiv k\!\!\pmod{3}$ for $1\leq j<\ell(\lambda)$}.
\end{align*}
\end{Conj}

Like the RR partition theorem, the Kanade-Russell conjecture has an equivalent $q$-series identity reformulation~\cite[\S3]{Kur}.
Because the double sums due to Kur\c{s}ung\"oz are well-known, we give a reformulation regarding triple sums.
Recall that 
Kanade-Russell also conjectured $\KR_4\PT T^{(9)}_{2,3,5,8}$ and
$\KR_5\PT T^{(9)}_{1,4,6,7}$~\cite{KR,Rus}, where
\begin{align*}
\KR_4 &= \{\lambda\in\PAR\mid \lambda_j-\lambda_{j+2}\geq 3 \textrm{ for $1\leq j\leq\ell(\lambda)-2$, $P(\lambda,2)$ holds and $m_1(\lambda)=0$}\},\\
\KR_5 &= \{\lambda\in\PAR\mid \lambda_j-\lambda_{j+2}\geq 3 \textrm{ for $1\leq j\leq\ell(\lambda)-2$, $P(\lambda,1)$ holds and $m_2(\lambda)\leq 1$}\}.
\end{align*}


\begin{Prop}\label{KRreformulation}
For $1\leq a\leq 5$, we have
\begin{align*}
\sum_{\lambda\in\KR_a}x^{\ell(\lambda)}q^{|\lambda|}=\sum_{i,j,k\geq 0}\frac{q^{3{i\choose 2}+8{j\choose 2}+6{k\choose 2}+4ij+3ik+6jk+A_{a,1}i+A_{a,2}j+A_{a,3}k}}{(q;q)_i(q^2;q^2)_j(q^3;q^3)_k}x^{i+2j+2k},
\end{align*}
where $A_1=(1,4,3)$, $A_{2}=(2,6,6)$, $A_{3}=(3,8,6)$, $A_{4}=(2,6,5)$ and $A_{5}=(1,6,4)$.
\end{Prop}
Because it is routine to derive a $q$-difference equation for the right (resp. left) hand side
by a $q$-analog of Sister Celine's technique~\cite{Rie} 
(resp. by Andrews' linked partition ideals~\cite{An0,An1} or the regularly linked sets by Takigiku and the author~\cite{TT}),
we omit an automatic proof of Proposition \ref{KRreformulation}.
See also \S\ref{automa}, ~\cite[\S7.1]{Tsu} and ~\cite[\S4]{CL}.

\subsection{The main result}

\begin{figure}
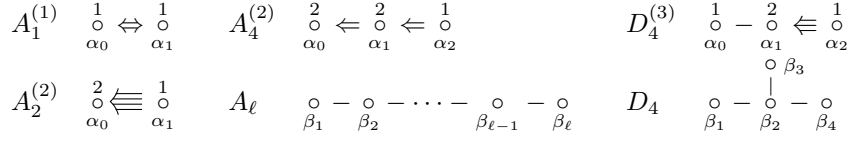

\[
\begin{array}{ll@{\qquad}ll@{\qquad}ll@{\qquad}ll}
A_1^{(1)} & \node{1}{\alpha_0} \Leftrightarrow \node{1}{\alpha_1} &
A_{4}^{(2)} &  \node{2}{\alpha_0}\Leftarrow \node{2}{\alpha_{1}}\Leftarrow\node{1}{\alpha_{2}} &
D_4^{(3)} & \node{1}{\alpha_0}-\node{2}{\alpha_1}\Lleftarrow\node{1}{\alpha_2} \\
A_2^{(2)} & \node{2}{\alpha_0} \Llleftarrow \node{1}{\alpha_1} &
A_{\ell} & \node{}{\beta_1} - \node{}{\beta_2} - \cdots - \node{}{\beta_{\ell-1}}-\node{}{\beta_\ell} &
D_4 & \node{}{\beta_1} - \node{\ver{}{\beta_{3}}}{\beta_{2}} -\node{}{\beta_{4}} 
\end{array}
\]
\caption{The Dynkin diagrams $A_1^{(1)}$, $A_2^{(2)}$, $A_4^{(2)}$, $D_4^{(3)}$ and $A_{\ell}$, $D_4$}
\label{twisted}
\end{figure}

Recall that for a standard module (a.k.a., integrable highest weight module) $V$ of an affine Kac-Moody Lie algebra~\cite{Kac},
the character $\chi(\Omega_V)$ of the vacuum space $\Omega_V$ with respect to the principal Heisenberg
subalgebra (see also \S\ref{zalg2}) is obtained by the Lepowsky numerator formula~\cite[Proposition 8.4]{LW3}. 

Let $i=1$ or $2$ and consider the level 3 standard module 
$V=V((i+1)\Lambda_0+(2-i)\Lambda_1)$ of the affine algebra $\GEE(A^{(1)}_1)$ (see Figure \ref{twisted}). 
Lepowsky-Milne observed that 
\begin{align*}
\chi(\Omega_V)=\frac{1}{(q^{i},q^{5-i};q^5)_{\infty}},
\end{align*}
is the infinite product in the RR identities~\cite{LM}.
In a seminal paper~\cite{LW3} (see also ~\cite{LW2}), Lepowsky-Wilson showed (see ~\cite[Theorem 10.4]{LW3})
\begin{align*}
\chi((\Omega_V)^{[n]}/(\Omega_V)^{[n-1]})=\frac{q^{n^2+n(i-1)}}{(q;q)_n},
\end{align*}
for $n\geq 0$, where $(\Omega_V)^{[n]}$ is the associated $Z$-filtration, which coincides with the $s$-filtration $(\Omega_V)_{[n]}$~\cite[\S5]{LW3}
(and $(\Omega_V)^{[-1]}=(\Omega_V)_{[-1]}=\{0\}$).
For that purpose, Lepowsky-Wilson proved that the set of $\Z$-monomials parameterized by $\RR_i$
\begin{align*}
\{\Z_{-\lambda_1}\cdots\Z_{-\lambda_{\ell}}v_0\mid (\lambda_1,\dots,\lambda_{\ell})\in\RR_i \}
\end{align*}
forms a basis of $\Omega_V$, where $v_0$ is a highest weight vector of $V$ and $\Z_i$ is the $Z$-operator (see ~\cite[(3.13)]{LW3} and a review in \S\ref{zalg}) associated with the root $\beta_1$ of $A_1$.

Our main result gives a weaker vertex operator interpretation of the Kanade-Russell conjecture
in the following way.

\begin{Thm}\label{mainres}
  For $1\leq a\leq 3$, the set of $\Z$-monomials parameterized by $\KR_a$
  \begin{align*}
    \{\Z_{-\lambda_1}\cdots\Z_{-\lambda_{\ell}}w_0\mid (\lambda_1,\dots,\lambda_{\ell})\in\KR_a \}
  \end{align*}
  spans $\Omega_{V(\Lambda^{(a)})}$, where $\Z_i$ is the $Z$-operator associated with the root $\beta_1$ of $D_4$ (see Figure \ref{twisted}),
  $w^{(a)}_0$ is a highest weight vector of the standard module $V(\Lambda^{(a)})$ of type $D^{(3)}_4$ whose
  highest weight is given by $\Lambda^{(1)}=\Lambda_0+\Lambda_1, \Lambda^{(2)}=3\Lambda_0$ and $\Lambda^{(3)}=\Lambda_2$.
\end{Thm}

An immediate consequence of Theorem \ref{mainres} and the infinite product expression of $\chi(\Omega_{V(\Lambda^{(a)})})$ is that the Kanade-Russell Conjecture (Conjecture \ref{krconj})
is equivalent to the claim that the $\Z$-monomials in Theorem \ref{mainres} are linearly independent. 

\begin{Cor}
We have $|\KR_a\cap\PAR(n)|\geq |T^{(9)}_{2^{a-1},3,6,9-2^{a-1}}\cap\PAR(n)|$ for $n\geq 0$ and $a=1,2,3$.
\end{Cor}

Our proof is a variant of 
~\cite[\S6]{LW3} and thus
standard in vertex operator theory (see also similar calculations in ~\cite{Cap0,Cap1,Ito,Kan,LW3,LW4,MP,Nan,Ta1,Ta2,Tsu,TX}) 
once
a suitable set of relations between the $Z$-operators is available.
In our case, four kinds of ``generalized commutation relations'' (see Theorem \ref{mainres2})
which are sometimes called with adjectives ``anti'' and ``partial'' (see ~\cite[\S3]{Hu2} and \S\ref{proofA}, \S\ref{proofB}, \S\ref{proofC}, \S\ref{proofD}) suffice.
It would be interesting to find a higher structure (e.g., vertex operator algebra structures as in ~\cite{Hu1,Hu2}) in our calculation
as well as to find a proof of linear independence (e.g., as in ~\cite{Cap1,LW3,MP,TX}).
We remark that
the arguments in \S\ref{forfour}, \S\ref{forfive}, \S\ref{forsix} are related to overlap ambiguities (resp. critical pairs)
in non-commutative Gr\"obner basis theory~\cite{Ber} (resp. theory of term rewriting~\cite{BN,New}).

\subsection{An $A^{(2)}_4$ analog}
The work of Lepowsky-Wilson~\cite{LW3}
initiated intensive research on 
explicit realizations of the generalized RR identity 
\begin{align*}
\chi(\Omega_V)=\sum_{n\geq 0}\chi((\Omega_V)^{[n]}/(\Omega_V)^{[n-1]})
\end{align*}
for a standard module $V$ of type $X^{(r)}_N=(a_{ij})_{i,j\in I}$ (see also ~\cite[Theorem 7.5]{LW3}). 
We give a brief review for types $A^{(1)}_1, A^{(2)}_2$ and $ A^{(1)}_2$ with levels greater than 2.


For $A^{(1)}_{1}$ arbitrary level, Lepowsky-Wilson showed that the $\Z$-monomials parameterized by partitions in
the Andrews-Gordon-Bressoud partition theorem (see ~\cite[\S3.2]{Sil} and the references therein) spans $\Omega_V$~\cite[Corollary 13.2, Lemma 14.3]{LW4}, which were proven to be linearly independent
by Meurman-Primc~\cite[\S9, Appendix]{MP}. 

For $A^{(2)}_{2}$ level 3 (resp. 4), Capparelli~\cite{Cap0} (resp. Nandi~\cite{Nan}) obtained a set of $\X$-monomials
that spans $V$ as a module over the principal Heisenberg algebra. Subsequently, they were proven to be linearly independent by ~\cite{An4,Cap1,TX} (resp. ~\cite{TT}).
On corresponding $q$-series identities, see ~\cite{BKRS,Ku2,TT2} and the references therein.

For $A^{(1)}_{2}$ level 3, the author showed that a set of $\X$-monomials parameterized by certain 2-color partitions
(or bipartitions) forms a basis of $V$ as a module over the principal Heisenberg algebra~\cite[Theorem 1.3]{Tsu}.

Although there is much literature devoted to higher levels, such
as ~\cite{Cap2,Hir,KR4,MS,TT2} (resp. ~\cite{ASW,CDA,CW,FFW,FW,Ka2,KR3,Tsu,Unc,War2,War,War3}) for $A^{(2)}_2$ (resp. $A^{(1)}_2$),
relations to representation theory remain unknown except the aforementioned levels. 
As in ~\cite[(6.8)]{LW3}, it is natural to expect that there exist certain $n(X^{(r)}_N)$-color partitions whose corresponding $\Z$-monomials form a basis of $\Omega_V$, where (see ~\cite[\S8, \S9]{Fig})
\begin{align*}
n(X^{(r)}_N)=\frac{\textrm{the number of roots of type $X_N$}}{\textrm{the $r$-twisted Coxeter number of $X_N$}}=|I|-1.
\end{align*}
For example, $n(A^{(1)}_1)=n(A^{(2)}_2)=1$ and $n(A^{(1)}_2)=n(A^{(2)}_4)=n(D^{(3)}_4)=2$.
Thus, the above results can be regarded as instances of that expectation.
On the other hand, for low levels such as level 2, sometimes $\Omega_V$ has a basis of $\Z$-monomials parameterized
by certain partitions as in ~\cite{Ito,Kan,Ta1,Ta2,Tsu,TX}.
While conjectural, Theorem \ref{mainres} may be thought of as a similar example. 

It would be interesting to find more examples similar to Theorem \ref{mainres} for other types $A$ with 
$n(A)>1$ (or other levels).
We note that, as shown in ~\cite{Ito},  ``the Kanade-Russell conjecture modulo 12''~\cite{KR2}, which is now a theorem~\cite{BJM,Ros}, 
can be understood as a kind of these examples for type $A^{(2)}_9$ level 2.

\begin{Thm}\label{mainconj}
Let $F$ be the set consisting of the following 13 partitions
\begin{align*}
(1,1,1),(2,1,1),(2,2,1),(3,2,1),(3,3,1),(5,3,3),(4,4,1,1),\\
(5,4,1,1),(5,4,2,1),(5,5,2,1),(6,5,3,1,1),(6,6,3,1,1),(7,6,4,2,1)
\end{align*}
and define three sets $I_1,I_2,I_3$ as follows.
\begin{align*}
I_1 &= \{(1),(5,4,2,2),(9,8,6,4,2,2)\},\quad
I_2 = \{(1,1),(2,2),(4,3,1)\},\\
I_3 &= \{(1,1),(2,1),(2,2),(3,2),(3,3),(4,3,1),(4,4,1),(5,4,2),(6,5,3,1)\}.
\end{align*}

For $1\leq a\leq 3$,
let $\EL_a$ be the set of partitions that does not begin with $\boldsymbol{c}$ for $\boldsymbol{c}\in I_a$
and does not match $(b_1+k,\dots,b_{p}+k)$ for $(b_1,\dots,b_{p})\in F$ and $k\geq 0$ (see \S\ref{automa}).
We have $\EL_1\PT T^{(16)}_{2,3,4,5,11,12,13,14}$,
$\EL_2\PT T^{(2)}_{1}$ and $\EL_3\PT T^{(16)}_{1,4,6,7,9,10,12,15}$.
\end{Thm}


We arrived at the statement of Theorem \ref{mainconj} by checking 
\begin{align*}
    \{\Z_{-\lambda_1}\cdots\Z_{-\lambda_{\ell}}u_0\mid (\lambda_1,\dots,\lambda_{\ell})\in\EL_a\cap\PAR(n) \}
\end{align*}
is a set of linearly independent vectors 
for $1\leq a\leq 3$ and $0\leq n\leq 17-2^{a-1}$.
Here, $\Z_i$ is the $Z$-operator associated with the root $\beta_1$ of $A_4$ (see Figure \ref{twisted}), and
$u_0$ is a highest weight vector of the standard module $V(\Upsilon^{(a)})$ of type $A^{(2)}_4$ whose
highest weight is given by $\Upsilon^{(1)}=3\Lambda_0, \Upsilon^{(2)}=\Lambda_0+\Lambda_1$ and $\Upsilon^{(3)}=\Lambda_0+\Lambda_2$.

After submission to arXiv of the first version of this paper, where Theorem \ref{mainconj} was stated as a conjecture,
we learned from Matthew Russell that Theorem \ref{mainconj} for $a=1$ is a theorem due to Andrews-van Ekeren-Heluani~\cite[Theorem 3]{AvH}.
In \S\ref{automaproof}, we point out that the case $a=2,3$ is also easily proven 
based on their work~\cite[Proposition 4.4]{AvH}. 
Interestingly, they seemed to find Theorem \ref{mainconj} for $a=1$ via the simple vertex algebra associated with the $(3,4)$ Virasoro minimal
model of central charge $c=1/2$, unlike via the level 3 standard modules of type $A^{(2)}_4$.


\hspace{0mm}

\noindent{\bf Organization of the paper.} 
We review the $Z$-operators in \S\ref{zalg}, 
establish the four relations among them in \S\ref{commrel} and
prove Theorem \ref{mainres} in \S\ref{mainprof}.
In \S\ref{automa}, 
we advertise the regularly linked sets~\cite{TT} by demonstrating an automatic derivation of
a $q$-difference equation for the 
generating function of $\EL_a$ (see Theorem \ref{mainqd}).
By using Theorem \ref{mainqd} and ~\cite[Proposition 4.4]{AvH}, we show Theorem \ref{mainconj} for $a=2,3$ in \S\ref{automaproof}.

\section{The $Z$-operators}\label{zalg}
In this section, we fix an affine Dynkin diagram $X^{(r)}_N$ of type $ADE$.

\subsection{The principal automorphisms}
Consider a lattice $L=\mathbb{Z}^N$ with a symmetric bilinear form $\langle,\rangle:L\times L\to \mathbb{Z}$ given by
$\langle\boldsymbol{x},\boldsymbol{y}\rangle={}^t\boldsymbol{x}X_N\boldsymbol{y}$.
Then, the standard basis $\{\boldsymbol{e}_1,\dots,\boldsymbol{e}_N\}$ of $L$, which we write as $\{\beta_1,\dots,\beta_N\}$,
is a set of simple roots for the set of roots $\Phi=\{\beta\in L\mid \langle\beta,\beta\rangle=2\}$.
The twisted Coxeter automorphism (see ~\cite[\S8]{Fig}) $\nu:L\to L$ is defined by $\nu=\sigma_{i_1}\cdots \sigma_{i_j}\sigma'$, where
\begin{enumerate}
\item $\sigma_i:L\to L,\boldsymbol{x}\mapsto \boldsymbol{x}-\langle \boldsymbol{x},\beta_i\rangle\beta_i$ is a reflection,
\item $\sigma':L\to L,\beta_{i}\mapsto\beta_{\sigma(i)}$ is the Dynkin diagram automorphism of order $r$
induced from that on the set of indices $\sigma:\{1,\dots,N\}\to\{1,\dots,N\}$, and
\item $\{i_1,\dots,i_j\}$ is a set of complete representatives for the orbits of the action
  of the cyclic group $\{\sigma^s\mid 0\leq s< r\}$ of order $r$ on the set of indices $\{1,\dots,N\}$.
\end{enumerate}

In the following, we fix an integer $m$ to the order of $\nu$ and fix a $m$-th primitive root of unity $\omega$.
As in ~\cite[(5.1)]{Fig}, we define a map $\varepsilon:L\times L\to\mathbb{C}^{\times}$ by
\begin{align*}
\varepsilon(\beta,\beta')=\prod_{p=1}^{m-1}(1-\omega^{-p})^{\langle\nu^p(\beta),\beta'\rangle}.
\end{align*}

Let $\GA=\mathbb{C}\otimes_{\mathbb{Z}}L$ be a complexification of $L$
and define a complex vector space $\GEE=\GA\oplus\bigoplus_{\beta\in\Phi}\mathbb{C}x_{\beta}$, where
$x_{\beta}$ is a formal symbol attached to a root $\beta$.
It is known that $\GEE$ with brackets determined by
\begin{align*}
[\beta_i,x_{\beta}] = \langle \beta_i,\beta\rangle x_{\beta} = -[x_{\beta},\beta_i],\quad
[x_{\beta},x_{\beta'}]=
\begin{cases}
  \varepsilon(-\beta,\beta)\beta & (\textrm{if $\langle\beta,\beta'\rangle=-2$}) \\
  \varepsilon(\beta,\beta')x_{\beta+\beta'} & (\textrm{if $\langle\beta,\beta'\rangle=-1$}) \\
  0 & (\textrm{if $\langle\beta,\beta'\rangle\geq 0$}),
\end{cases}
\end{align*}
where $\beta,\beta'\in\Phi$ and $1\leq i\leq N$,
is isomorphic to the Kac-Moody Lie algebra $\GEE(X_N)$ of type $X_N$.
The form $\langle,\rangle$ is extended to $\GEE$ by
$\langle \beta_i,x_{\beta}\rangle=0=\langle x_{\beta},\beta_i\rangle$ and 
$\langle x_{\beta},x_{\beta'}\rangle=\varepsilon(\beta,\beta')\delta_{\beta+\beta',0}$.
The map $\nu$ is also extended to $\GEE$ by $\nu(x_{\beta})=x_{\nu(\beta)}$.
Then, $\nu$ is a principal automorphism of $\GEE$, and
$\langle,\rangle$ is a non-degenerate invariant form which is $\nu$-invariant (see ~\cite[\S6, \S9]{Fig}).

\subsection{The Lepowsky-Wilson $Z$-algebras}\label{zalg2}
The $\nu$-twisted affinization $\GE$ given by 
\begin{align*}
  \GE=\bigoplus_{i\in\mathbb{Z}}(\GEE_{(i)}\otimes t^{i})\oplus\mathbb{C}c\oplus\mathbb{C}d,
\end{align*}
where $\GEE_{(i)}=\{x\in\GEE\mid \nu^i(x)=\omega^i x\}$, with brackets determined by
\begin{align*}
  [x\otimes t^i,y\otimes t^j]=[x,y]\otimes t^{i+j}+\frac{i\delta_{i+j,0}}{m}\langle x,y\rangle c,\quad
  [d,x\otimes t^i]=tx\otimes t^i,
\end{align*}
and by the condition that $c$ is central, is isomorphic to $\GEE(X^{(r)}_N)$ (see ~\cite[Chapter 8]{Kac} and ~\cite[Proposition 9.4]{Fig}).

For $k\in\mathbb{C}^{\times}$, as in ~\cite[\S3]{LW3}, let $\LWC_k$ be the full subcategory of the category of $\GE$-modules whose objects $V$ satisfy
the following three conditions.
\begin{enumerate}
\item $c$ acts as a scalar multiplication by $k$, i.e., $V$ has level $k$.
\item $V$ has a simultaneous eigenspace decomposition $V=\bigoplus_{\lambda\in\GH^{\ast}}V_{\lambda}$ with respect to the Cartan subalgebra $\GH=(\GEE_{(0)}\otimes 1)\oplus\mathbb{C}c\oplus\mathbb{C}d$ (see ~\cite[\S9.3]{Fig}).
\item For $z\in\mathbb{C}$, there exists $i_0\in\mathbb{Z}$ such that $V_{z+i}=\{0\}$ for $i>i_0$.
\end{enumerate}

For an object $V$ in $\LWC_k$, the Lepowsky-Wilson $Z$-algebra $Z_V$ (in the principal picture)
is defined to be the subalgebra of $\END V$ generated by $c$, $d$ and $Z_i(\beta)$
for $i\in\mathbb{Z}$ and $\beta\in\Phi$~\cite[p.222]{LW3}, where $\PR_i:\GEE\to\GEE_{(i)}$ for $i\in\mathbb{Z}$ is the projection and
\begin{align}
Z(\beta,\zeta) &= \sum_{i\in\mathbb{Z}}Z_i(\beta)\zeta^{i}  = E^{-}(\beta,\zeta,k)X(\beta,\zeta)E^{+}(\beta,\zeta,k),\\
X(\beta,\zeta) &= \sum_{i\in\mathbb{Z}}(\PR_i(x_{\beta})\otimes t^i)\zeta^i,\\
E^{\pm}(\beta,\zeta,r) &= \sum_{\pm i\geq 0}E_i^{\pm}(\beta)\zeta^i
=\exp\left(m\sum_{\pm j>0}\frac{\PR_j(\beta)\otimes t^j}{rj}\zeta^j\right).
\label{basi2}
\end{align}

As in ~\cite[Proposition 4.7]{LW3}, the vacuum space $\Omega_V$ is defined as
\begin{align*}
\Omega_V=\{v\in V\mid \textrm{$av=0$ for $a\in\GAAA_{+}$}\},
\end{align*}
where $\GAAA=[\GAA,\GAA]= \GAAA_{+}\oplus\GAAA_{-}\oplus\mathbb{C}c$ is the principal Heisenberg subalgebra
defined by
\begin{align*}
\GAA = \GAAA_{+}\oplus\GAAA_{-}\oplus\mathbb{C}c\oplus\mathbb{C}d,\quad
\GAAA_{\pm}=\bigoplus_{\pm i>0}(\GEE_{(i)}\cap \GA)\otimes t^i.
\end{align*}


For a complex vector space $U$, let $U\{\zeta_1,\dots,\zeta_r\}$ be the complex vector space of formal Laurent series
in the variables $\zeta_1,\dots,\zeta_r$ with coefficients in $U$, i.e.,
\begin{align*}
  \textstyle
  U\{\zeta_1,\dots,\zeta_r\} = \{
  \sum_{n_1,\dots,n_r\in\mathbb{Z}}u_{n_1,\dots,n_r}\zeta_1^{n_1}\cdots\zeta_r^{u_r}\mid u_{n_1,\dots,n_r}\in U\}.
\end{align*}
The equalities in the following citations are those in $\END(V)\{\zeta_1,\zeta_2\}$ or $\END(V)\{\zeta\}$,
where $V$ is an object in $\LWC_k$, $\beta,\beta'\in\Phi$ and $r,s\geq 1$, $p\in\mathbb{Z}$.
As usual, the formal delta is defined by $\delta(\zeta)=\sum_{n\in\mathbb{Z}}\zeta^n$ and let $D=\zeta\cdot d/d\zeta$
so that $D\delta(\zeta)=\sum_{n\in\mathbb{Z}}n\zeta^n$.

\begin{Thm}[{\cite[Theorem 3.10]{LW3}, \cite[Theorem 7.3]{Fig}}]\label{gencomm}
\begin{align*}
{} &F_{\beta,\beta'}(\zeta_1/\zeta_2)Z(\beta,\zeta_1)Z(\beta',\zeta_2)-F_{\beta',\beta}(\zeta_2/\zeta_1)Z(\beta',\zeta_2)Z(\beta,\zeta_1)\\
&=
\frac{1}{m}\sum_{p\in C_{-1}}\varepsilon(\nu^p(\beta),\beta')Z(\nu^p(\beta)+\beta',\zeta_2)\delta(\omega^{-p}\zeta_1/\zeta_2)
+\frac{k\langle x_{\beta},x_{-\beta}\rangle}{m^2}\sum_{p\in C_{-2}}D\delta(\omega^{-p}\zeta_1/\zeta_2)
\end{align*}
where $C_{\ell}=\{0\leq p<m\mid \langle \nu^p(\beta),\beta'\rangle=\ell\}$ and
$F_{\beta,\beta'}(x)=\prod_{p=0}^{m-1}(1-\omega^{-p}x)^{\langle \nu^p(\beta),\beta'\rangle/k}$.
\end{Thm}

\begin{Prop}[{\cite[Proposition 3.4]{LW3}}]\label{gencomm1}
\begin{align*}
E^{+}(\beta,\zeta_1,r)E^{-}(\beta',\zeta_2,s)
=E^{-}(\beta',\zeta_2,s)E^{+}(\beta,\zeta_1,r)F_{\beta,\beta'}(\zeta_1/\zeta_2)^{1/rs}.
\end{align*}
\end{Prop}

\begin{Prop}[{\cite[Theorem 3.3]{LW3}, \cite[Theorem 7.2]{Fig}}]\label{gencomm2}
$Z(\nu^p(\beta),\zeta)=Z(\beta,\omega^p\zeta)$.
\end{Prop}

\begin{Prop}[{\cite[Proposition 3.2.(3)]{LW3}}]\label{gencomm3}
  $DE^{\pm}(\beta,\zeta,r)=\frac{m}{r}\beta^{\pm}(\zeta)E^{\pm}(\beta,\zeta,r)$,
  where $\beta^{\pm}(\zeta)=\sum_{\pm j>0}(\PR_j(\beta)\otimes t^j)\zeta^j$.
\end{Prop}


\subsection{The basic representations}\label{basrep}

As shown in ~\cite[Theorem 9.7]{Fig}, giving (an isomorphism class of) a basic representation of $\GE$
is the same as giving
a $\nu$-invariant (i.e., $\rho(\nu(\beta))=\rho(\beta)$ for $\beta\in L$) group homomorphism $\rho:L\to\mathbb{C}^{\times}$.
The corresponding basic representation $V^{(\rho)}$ is given
by the underlying $\GAA$-module $V=\IND_{\GAAA_{+}\oplus\mathbb{C}c\oplus\mathbb{C}d}^{\GAA}\mathbb{C}(\cong U(\GAAA_{-}))$, where 
$\GAAA_{+}\oplus\mathbb{C}d$ (resp. $c$) acts
trivially (resp. as the identity) on $\mathbb{C}$ and the action by
$X(\beta,\zeta)=\frac{\rho(\beta)}{m}E^{-}(-\beta,\zeta,1)E^{+}(-\beta,\zeta,1)$.


\section{The (partial) generalized (anti)commutation relations}\label{commrel}
\subsection{The affine type $D^{(3)}_4$}\label{dyn}
We apply the general theory reviewed in \S\ref{zalg} to $D^{(3)}_4$ (see Figure \ref{twisted}).
Our convection of the Dynkin diagram automorphism $\sigma$ is given by
$\sigma(1)=3,\sigma(3)=4,\sigma(4)=1,\sigma(2)=2$,
and that of the twisted Coxeter automorphism is given by $\nu=\sigma_1\sigma_2\sigma'$.

Then, the order $m$ of $\nu$ is $12$, and the set of roots $\{\beta^{(i)}_{j}\mid \textrm{$i=1,2$ and $0\leq j< m$}\}$,
where $\beta^{(i)}_{j}=\nu^j(\beta_i)$ is given by $\nu^6=-\ID$ and
\begin{align*}
\begin{array}{lll}
\beta^{(1)}_1=\beta_1+\beta_2+\beta_3, &
\beta^{(1)}_2=\beta_1+\beta_2+\beta_3+\beta_4, &
\beta^{(1)}_3=\beta_1+2\beta_2+\beta_3+\beta_4, \\
\beta^{(1)}_4=\beta_2+\beta_3+\beta_4, &
\beta^{(1)}_5=\beta_2+\beta_4, &
\beta^{(2)}_1=-\beta_1-\beta_2, \\
\beta^{(2)}_2=-\beta_3,\quad \beta^{(2)}_5=-\beta_4, &
\beta^{(2)}_3=-\beta_1-\beta_2-\beta_4, &
\beta^{(2)}_4=-\beta_2-\beta_3.
\end{array}
\end{align*}

\subsection{Summary}

Recall $m=12$. 
As in \S\ref{zalg}, $\omega$ is 
a fixed primitive $m$-th root of unity, a solution of the $m$-th cyclotomic polynomial $\Phi_{m}(x)=x^4-x^2+1$.
Let
\begin{align*}
\AD &= 4+4\omega-2\omega^3,\quad
\BD = -24-28\omega+14\omega^3,\quad
\CD = 42+48\omega-24\omega^3,\\
\A &= -6-8\omega+4\omega^3,\quad
\B = 14+16\omega-8\omega^3,\quad
\D = 4-8\omega^2-6\omega^3,\quad
\E = 2-4\omega^2. 
\end{align*}

\begin{Def}
For integers $a_0,\dots,a_5$, we define a formal power series 
\begin{align*}
H_{a_0,\dots,a_5}(x)=\prod_{p=0}^{5}(1-\omega^{-p}x)^{a_p/3}(1+\omega^{-p}x)^{-a_p/3}.
\end{align*}
We also define $G_i(x)=\sum_{p\geq 0}c^{(i)}_px^p$ for $1\leq i\leq 6$ by
\begin{align*}
  G_1 &= H_{2,1,1,0,-1,-1},\quad
  G_2 = H_{-1,1,1,0,-1,-1},\quad
  G_3 = H_{-1,1,-2,0,2,-1},\\ 
  G_4 &= H_{2,-2,-2,0,2,2},\quad
  G_5 = H_{2,-2,1,0,-1,2},\quad
  G_6 = G_4-\frac{\D}{\E}G_5.
\end{align*}
\end{Def}

The four relations in Theorem \ref{mainres2} will be proven in this order in \S\ref{proofA}, \S\ref{proofB}, \S\ref{proofC}, \S\ref{proofD} 
after preparation in \S\ref{compan}. 
\begin{Thm}\label{mainres2}
Let $W=(V^{(\rho)})^{\otimes 3}$ be the triple tensor product of the basic representation $V^{(\rho)}$ of type $D^{(3)}_4$,
where $\rho:\mathbb{Z}^4\to\mathbb{C}^{\times}$ is the trivial group homomorphism.
Put $\Z_i=Z_i(\beta_1)$ and $\Z'_i=Z_i(\beta_2)$.
In the $Z$-algebra $Z_W$, for $A,B\in\mathbb{Z}$, we have 
\begin{align*}
{} &\sum_{p\geq 0}c^{(1)}_p(\Z_{A-p}\Z_{B+p}-\Z_{B-p}\Z_{A+p}) \\
&=
\frac{\varepsilon(\nu^4\beta_1,\beta_1)}{m}(\omega^{-2A+2B}-\omega^{2A-2B})\Z_{A+B}
+\frac{\varepsilon(\nu^5\beta_1,\beta_1)}{m}(\omega^{4A+9B}-\omega^{9A+4B})\Z'_{A+B}+\frac{3\varepsilon(\beta_1,-\beta_1)}{m^2}\delta_{A+B,0}A(-1)^A,\\
{} &\sum_{p\geq 0}c^{(2)}_p(\Z_{A-p}\Z_{B+p}+\Z_{B-p}\Z_{A+p}) \\
&=
\frac{\AD}{m}(\omega^{-2A+2B}+\omega^{2A-2B})\Z_{A+B}
+\frac{\BD}{m}(\omega^{4A+9B}+\omega^{9A+4B})\Z'_{A+B}
+\frac{3\CD}{m^2}\delta_{A+B,0}(-1)^A + \frac{4}{m}(-1)^{A+B}\Z_{A+B}.
\end{align*}

If $A+B$ is not divisible by 3, then we have
\begin{align*}
\sum_{p\geq 0}c^{(3)}_p(\Z_{A-p}\Z_{B+p}+\Z_{B-p}\Z_{A+p}) 
&=
\frac{\A}{m}(\omega^{4A+9B}+\omega^{9A+4B})\Z'_{A+B}
+\frac{3\B}{m^2}\delta_{A+B,0}(-1)^A + \frac{12}{m}(-1)^{A+B}\Z_{A+B},\\
\sum_{p\geq 0}c^{(6)}_p(\Z_{A-p}\Z_{B+p}-\Z_{B-p}\Z_{A+p}) 
&=
\frac{\D}{m}(\omega^{-2A+2B}-\omega^{2A-2B})\Z_{A+B}
+\frac{12}{m^2}(1-\frac{3\D}{\E})\delta_{A+B,0}(-1)^A A.
\end{align*}
\end{Thm}


\subsection{Auxiliary Fourier expansions}\label{compan}
Let $\BARR{\phantom{a}}:\mathbb{C}\{x\}\to \mathbb{C}\{x\}$ be a $\mathbb{C}$-linear map defined by $\BARR{\sum_{p\in\mathbb{Z}}c_px^p}=\sum_{p\in\mathbb{Z}}c_px^{-p}$,
which is also a $\mathbb{C}$-algebra homomorphism when restricted to $\mathbb{C}[[x]]$ (or $\mathbb{C}[[x^{-1}]]$).
The following results are elementary complex analysis exercises, and we omit the detail (see also ~\cite[Lemma 3.2, Lemma 3.5]{TX}).

\begin{Prop}\label{compan2}
\begin{align*}
G_1^{2}G_2+\BARR{G_1^{2}G_2} &= \AD(\delta(\omega^4 x)+\delta(\omega^{-4}x))+\BD(\delta(\omega^5 x)+\delta(\omega^{-5}x))+\CD\delta(-x),\\
G_1^{-1}G_2+\BARR{G_1^{-1}G_2} &= 2\delta(x),\\
G_1^{2}G_3+\BARR{G_1^{2}G_3} &= \A(\delta(\omega^5 x)+\delta(\omega^{-5}x))+\B\delta(-x),\\
G_1^{-1}G_3+\BARR{G_1^{-1}G_3} &= 6\delta(x)-2(\delta(\omega^2 x)+\delta(\omega^{-2} x)),\\
G_1^2G_4-\BARR{G_1^2G_4} &= 4D\delta(-x),\\
G_1^{-1}G_4-\BARR{G_1^{-1}G_4} &= \D(\delta(\omega x)-\frac{1}{3}\delta(\omega^2 x)+\frac{1}{3}\delta(\omega^{-2}x)-\delta(\omega^{-1}x)),\\
G_1^{2}G_5-\BARR{G_1^{2}G_5} &= 12D\delta(-x)+\E(\delta(\omega^4 x)-\delta(\omega^{-4}x)),\\
G_1^{-1}G_5-\BARR{G_1^{-1}G_5} &= \E(\delta(\omega x)-\delta(\omega^{-1}x)).
\end{align*}
\end{Prop}

\subsection{A proof of the generalized commutation relation}\label{proofA}
Apply Theorem \ref{gencomm} for $k=3$ and $\beta=\beta'=\beta_1$ by using $F_{\beta,\beta'}=G_1$ and $C_{-1}=\{4,5,7,8\}$ with
\begin{align*}
\nu^{\pm 4}(\beta_1)+\beta_1 &= \nu^{\pm 2}(\beta_1),\quad
\nu^{5}(\beta_1)+\beta_1=\nu^{9}(\beta_2),\quad
\nu^{7}(\beta_1)+\beta_1=\nu^{4}(\beta_2),\\
\varepsilon(\nu^4(\beta_1),\beta_1) &= \D=-\varepsilon(\nu^{-4}(\beta_1),\beta_1),\quad
\varepsilon(\nu^5(\beta_1),\beta_1) = -52+104\omega^2+90\omega^3=-\varepsilon(\nu^7(\beta_1),\beta_1).
\end{align*}
We get the result
by comparing the coefficients of $\zeta_1^A\zeta_2^B$ of the identity in Theorem \ref{gencomm} using Proposition \ref{gencomm2}.

\subsection{A proof of the generalized anticommutation relation}\label{proofB}
As in ~\cite[Lemma 9.1]{LW3}, for $\gamma\in\Phi$, we have
$Z(\gamma,\zeta)=Z^{(1)}(\gamma,\zeta)+Z^{(2)}(\gamma,\zeta)+Z^{(3)}(\gamma,\zeta)$,
where, for $1\leq j\leq 3$, $mZ^{(j)}(\gamma,\zeta)$ is the tensor product of three factors, the $i(\ne j)$-th (resp. $j$-th) tensorand being (see also ~\cite[(3.6), (3.7), (3.8)]{TX} and \S\ref{zalg2}, \S\ref{basrep})
\begin{align*}
E^{-}(\gamma,\zeta,3)E^{+}(\gamma,\zeta,3)\quad \textrm{ (resp. $E^{-}(\gamma,\zeta,3)^{-2}E^{+}(\gamma,\zeta,3)^{-2}$)}.
\end{align*}

Recall $F_{\beta_1,\beta_1}=G_1$ as in \S\ref{proofA}, where $F_{\beta,\beta'}$ is defined as in Theorem \ref{gencomm} for $k=3$.
By Proposition \ref{gencomm1} together with 
\begin{align*}
\frac{1}{3}\cdot\frac{1}{3}\cdot2+\frac{-2}{3}\cdot\frac{-2}{3}=\frac{2}{3},\quad
\frac{1}{3}\cdot\frac{-2}{3}+\frac{-2}{3}\cdot\frac{1}{3}+\frac{1}{3}\cdot\frac{1}{3}=-\frac{1}{3},
\end{align*}
we see the following two results.

First, as in ~\cite[(3.26)]{TX}, we have $Z^{(j)}(\beta_1,\zeta_1)Z^{(j)}(\beta_1,\zeta_2)=\frac{1}{m^2}F_{\beta_1,\beta_1}(\zeta_1/\zeta_2)^2\QUE^{(j)}$, where $\QUE^{(j)}$ is the tensor product of three factors,
the $i$-th tensorand being
\begin{align*}
  E^{-}(\beta_1,\zeta_1,3)^{1-3\delta_{ij}}E^{-}(\beta_1,\zeta_2,3)^{1-3\delta_{ij}}
  E^{+}(\beta_1,\zeta_1,3)^{1-3\delta_{ij}}E^{+}(\beta_1,\zeta_2,3)^{1-3\delta_{ij}}.
\end{align*}
Here, $1\leq i\leq 3$ and $\delta_{ij}$ is the Kronecker delta.

Second, for $1\leq s\ne t\leq 3$, take a unique $1\leq u\leq 3$ such that $\{s,t,u\}=\{1,2,3\}$. Then,
as in ~\cite[(3.34)]{TX}, we have $Z^{(s)}(\beta_1,\zeta_1)Z^{(t)}(\beta_1,\zeta_2)=\frac{1}{m^2}F_{\beta_1,\beta_1}(\zeta_1/\zeta_2)^{-1}\QUEE^{(s,t)}$, where $\QUEE^{(s,t)}$ is the tensor product of three factors,
the $s$-th tensorand being
\begin{align*}
E^{-}(\beta_1,\zeta_1,3)^{-2}E^{-}(\beta_1,\zeta_2,3)E^{+}(\beta_1,\zeta_1,3)^{-2}E^{+}(\beta_1,\zeta_2,3),
\end{align*}
the $t$-th tensorand being
\begin{align*}
E^{-}(\beta_1,\zeta_1,3)E^{-}(\beta_1,\zeta_2,3)^{-2}E^{+}(\beta_1,\zeta_1,3)E^{+}(\beta_1,\zeta_2,3)^{-2},
\end{align*}
and the $u$-th tensorand being
\begin{align*}
E^{-}(\beta_1,\zeta_1,3)E^{-}(\beta_1,\zeta_2,3)E^{+}(\beta_1,\zeta_1,3)E^{+}(\beta_1,\zeta_2,3).
\end{align*}

Let $\QUE=\QUE^{(1)}+\QUE^{(2)}+\QUE^{(3)}$ and $\QUEE=\sum_{1\leq s\ne t\leq 3}\QUEE^{(s,t)}$.
By Proposition \ref{compan2} and in virtue of ``the residue theorem''~\cite[Proposition 3.9.(2)]{LW3},
we see as in ~\cite[(3.26), (3.34)]{TX}
\begin{align*}
{} &{}  G_2(\zeta_1/\zeta_2)Z(\beta_1,\zeta_1)Z(\beta_1,\zeta_2)+G_2(\zeta_2/\zeta_1)Z(\beta_1,\zeta_2)Z(\beta_1,\zeta_1)\\
&=  \frac{\AD}{m^2}(\delta(\omega^4 \zeta_1/\zeta_2)\QUE(\omega^{-4}\zeta_2,\zeta_2)+\delta(\omega^{-4} \zeta_1/\zeta_2)\QUE(\omega^{4}\zeta_2,\zeta_2))
  +\frac{\CD}{m^2}\delta(-\zeta_1/\zeta_2)\QUE(-\zeta_2,\zeta_2)\\
&\quad  +\frac{\BD}{m^2}(\delta(\omega^5 \zeta_1/\zeta_2)\QUE(\omega^{-5}\zeta_2,\zeta_2)+\delta(\omega^{-5}\zeta_1/\zeta_2)\QUE(\omega^{5}\zeta_2,\zeta_2))
+\frac{2}{m^2}\delta(\zeta_1/\zeta_2)\QUEE(\zeta_2,\zeta_2).
\end{align*}

It is easy to see that $\QUE^{(j)}(-\zeta_2,\zeta_2)=1$ (see also ~\cite[(3.29)]{TX})
and $\QUEE^{(s,t)}(\zeta_2,\zeta_2)=mZ^{(u)}(-\beta_1,\zeta_2)$ (see also ~\cite[(3.36)]{TX}).
It is also easy to see that
\begin{align*}
\QUE^{(j)}(\omega^{-4}\zeta_2,\zeta_2) &= mZ^{(j)}(\nu^{-2}(\beta_1),\zeta_2),\quad
\QUE^{(j)}(\omega^{4}\zeta_2,\zeta_2)=mZ^{(j)}(\nu^2(\beta_1),\zeta_2),\\
\QUE^{(j)}(\omega^{-5}\zeta_2,\zeta_2) &= mZ^{(j)}(\nu^4(\beta_2),\zeta_2),\quad
\QUE^{(j)}(\omega^{5}\zeta_2,\zeta_2)=mZ^{(j)}(\nu^9(\beta_2),\zeta_2),
\end{align*}
by Proposition \ref{gencomm2} and $\nu^{\mp 4}(\beta_1)+\beta_1=\nu^{\mp 2}(\beta_1)$, 
$\nu^{-5}(\beta_1)+\beta_1=\nu^4(\beta_2)$, $\nu^{5}(\beta_1)+\beta_1=\nu^9(\beta_2)$ (see also ~\cite[(3.27), (3.28)]{TX} and \S\ref{dyn}).
Thus, we have
\begin{align*}
{} &{} G_2(\zeta_1/\zeta_2)Z(\beta_1,\zeta_1)Z(\beta_1,\zeta_2)+G_2(\zeta_2/\zeta_1)Z(\beta_1,\zeta_2)Z(\beta_1,\zeta_1) \\
&= \frac{\AD}{m}\delta(\omega^{-4}\zeta_1/\zeta_2)Z(\nu^{2}(\beta_1),\zeta_2) + \frac{\AD}{m}\delta(\omega^{4}\zeta_1/\zeta_2)Z(\nu^{-2}(\beta_1),\zeta_2)+\frac{3\CD}{m^2}\delta(-\zeta_1/\zeta_2) \\
&+ \frac{\BD}{m}\delta(\omega^{-5}\zeta_1/\zeta_2)Z(\nu^{9}(\beta_2),\zeta_2) +\frac{\BD}{m}\delta(\omega^{5}\zeta_1/\zeta_2)Z(\nu^{4}(\beta_2),\zeta_2)+\frac{4}{m}\delta(\zeta_1/\zeta_2)Z(-\beta_1,\zeta_2).
\end{align*}
We get the result by comparing the coefficients of $\zeta_1^A\zeta_2^B$ as in \S\ref{proofA}.

\subsection{A proof of the partial generalized anticommutation relation}\label{proofC}
As in \S\ref{proofB}, $G_3(\zeta_1/\zeta_2)Z(\beta_1,\zeta_1)Z(\beta_1,\zeta_2)+G_3(\zeta_2/\zeta_1)Z(\beta_1,\zeta_2)Z(\beta_1,\zeta_1)$
is expanded as
\begin{align*}
{} &{} \frac{3\B}{m^2}\delta(-\zeta_1/\zeta_2)+
\frac{\A}{m}(\delta(\omega^{-5}\zeta_1/\zeta_2)Z(\nu^{9}(\beta_2),\zeta_2) +
\delta(\omega^{5}\zeta_1/\zeta_2)Z(\nu^{4}(\beta_2),\zeta_2))\\
&+ \frac{12}{m}\delta(\zeta_1/\zeta_2)Z(-\beta_1,\zeta_2)
-\frac{2}{m^2}(\delta(\omega^2 \zeta_1/\zeta_2)\QUEE(\omega^{-2}\zeta_2,\zeta_2)+\delta(\omega^{-2} \zeta_1/\zeta_2)\QUEE(\omega^{2}\zeta_2,\zeta_2)).
\end{align*}

Let $\QUEE(\omega^{\pm 2}\zeta_2,\zeta_2)=\sum_{i\in\mathbb{Z}}\QUEEE^{\pm}_{i}\zeta_2^i$.
In order to complete the proof, it is enough to show $\QUEEE^{\pm}_i=0$ when $i\not\in 3\mathbb{Z}$. The argument is the same as ~\cite[Corollary 3.9]{TX}. We duplicate it below for completeness.

We see that
$\QUEEE^{(s,t)}(\zeta_2)=\sum_{i\in\mathbb{Z}}\QUEEE^{(s,t)}_i\zeta_2^i$ defined as $\QUEE^{(s,t)}(\omega^{-2}\zeta_2,\zeta_2)$ is the tensor product of three factors,
the $v$-th tensorand being $E^{-}(\gamma_v,\zeta_2,3)E^{+}(\gamma_v,\zeta_2,3)$, where
\begin{align*}
\gamma_s &= -2\nu^{-2}(\beta_1)+\beta_1=\beta_1+2\beta_2+2\beta_3+2\beta_4,\\
\gamma_t &= \nu^{-2}(\beta_1)-2\beta_1=-2\beta_1-\beta_2-\beta_3-\beta_4,\\
\gamma_u &= \nu^{-2}(\beta_1)+\beta_1=\beta_1-\beta_2-\beta_3-\beta_4
\end{align*}
for $v=s,t,u$.
Note that $\gamma_t=\nu^4(\gamma_s)$ and $\gamma_u=\nu^{-4}(\gamma_s)$.
Thus, we have 
\begin{align*}
\QUEEE^{(3,1)}(\zeta_2)=\QUEEE^{(1,2)}(\omega^4\zeta_2),\quad
\QUEEE^{(2,3)}(\zeta_2)=\QUEEE^{(1,2)}(\omega^{-4}\zeta_2),\\
\QUEEE^{(3,2)}(\zeta_2)=\QUEEE^{(2,1)}(\omega^4\zeta_2),\quad
\QUEEE^{(1,3)}(\zeta_2)=\QUEEE^{(2,1)}(\omega^{-4}\zeta_2),
\end{align*}
and $\QUEE(\omega^{-2}\zeta_2,\zeta_2)=\sum_{i\in\mathbb{Z}}(1+\omega^{4i}+\omega^{-4i})(\QUEEE_i^{(1,2)}+\QUEEE_i^{(2,1)})\zeta_2^i$, which shows
$\QUEEE^{-}_i=0$ when $i\not\in 3\mathbb{Z}$.
The proof of $\QUEEE^{+}_i=0$ when $i\not\in 3\mathbb{Z}$ is the same.


\subsection{A proof of the partial generalized commutation relation}\label{proofD}
We note $D\delta(-\zeta_1/\zeta_2)\QUE^{(j)}=D\delta(-\zeta_1/\zeta_2)$ by the same argument as ~\cite[Proposition 3.10]{TX}
in virtue of ``the residue theorem''~\cite[Proposition 3.9.(3)]{LW3}, Proposition \ref{gencomm3}
and $\frac{-2}{3}+\frac{1}{3}+\frac{1}{3}=0$ (see also ~\cite[Lemma 3.6]{TX}).
Thus, as in \S\ref{proofC}, we have 
\begin{align*}
{} &{} G_6(\zeta_1/\zeta_2)Z(\beta_1,\zeta_1)Z(\beta_1,\zeta_2)-G_6(\zeta_2/\zeta_1)Z(\beta_1,\zeta_2)Z(\beta_1,\zeta_1)\\
&= \frac{12}{m^2}(1-\frac{3\D}{\E})D\delta(-\zeta_1/\zeta_2)+
\frac{\D}{m}(\delta(\omega^{-4}\zeta_1/\zeta_2)Z(\nu^{2}(\beta_1),\zeta_2) -\delta(\omega^{4}\zeta_1/\zeta_2)Z(\nu^{-2}(\beta_1),\zeta_2))\\
&\quad+ \frac{\D}{3m^2}(\delta(\omega^{-2} \zeta_1/\zeta_2)\QUEE(\omega^{2}\zeta_2,\zeta_2)-\delta(\omega^{2} \zeta_1/\zeta_2)\QUEE(\omega^{-2}\zeta_2,\zeta_2)),
\end{align*}
because by Proposition \ref{compan2} we know
\begin{align*}
  G_1^2G_6 - \BARR{G_1^2G_6} &= 4(1-\frac{3\D}{\E})D\delta(-x)-\D(\delta(\omega^4 x)-\delta(\omega^{-4} x)),\\
  G_1^{-1}G_6 - \BARR{G_1^{-1}G_6} &= \frac{\D}{3}(\delta(\omega^{-2} x)-\delta(\omega^{2} x)).
\end{align*}
We get the result by comparing the coefficients of $\zeta_1^A\zeta_2^B$ as in \S\ref{proofC}.

\section{A proof of Theorem \ref{mainres}}\label{mainprof}
Recall $\Lambda^{(a)}$ and $w^{(a)}_0$ for $1\leq a\leq 3$ in Theorem \ref{mainres} and 
let $V=V(\Lambda^{(a)})$. By $\chi(\Omega_{V(\Lambda^{(a)})})=1/(q^{2^{a-1}},q^3,q^6,q^{9-2^{a-1}};q^9)_{\infty}$, we have $\Z_{-1}w^{(2)}_0=0$, $\Z_{-1}w^{(3)}_0=0$ and $\Z_{-2}w^{(3)}_0=0$.

\subsection{An elimination of $\Z'_i$}
In virtue of ~\cite[Theorem 7.1]{LW3}, we have $\Omega_V=Z_{V}w^{(a)}_0$.
By substituting $A=B$ or $A=B+1$ of the second equality in Theorem \ref{mainres2}
(i.e., the generalized anticommutation relation proven in \S\ref{proofB}), we see that
$\Omega_V=Z^{(1)}_{V}w^{(a)}_0$ holds, where $Z^{(1)}_V$ is a subalgebra of $Z_V$ generated by $\Z_i$ for $i\in\mathbb{Z}$. 

\subsection{A Synopsis}

\begin{Def}[{\cite[\S6]{LW3}, \cite[\S4]{TX}}]
For $\boldsymbol{i}=(i_1,\dots,i_{\ell})\in\mathbb{Z}^{\ell}$ and $\boldsymbol{j}=(j_1,\dots,j_{\ell'})\in\mathbb{Z}^{\ell'}$, where $\ell$ and $\ell'$ are
non-negative integers,
we say that $\boldsymbol{i}$ is higher than $\boldsymbol{j}$ (denoted by $\boldsymbol{i}>\boldsymbol{j}$)
if $\ell<\ell'$ or 
\begin{align*}
  \ell=\ell', \boldsymbol{i}\ne\boldsymbol{j} \textrm{ and }
  i_{p}+\cdots+i_{\ell}\geq j_{p}+\cdots+j_{\ell} \textrm{ for $1\leq p\leq \ell$}.
\end{align*}
\end{Def}

For $\boldsymbol{j}=(j_1,\dots,j_{\ell'})\in\mathbb{Z}^{\ell'}$, 
we define $T(\boldsymbol{j})$ (resp. $T^{=}(\boldsymbol{j})$), which is also denoted by $T(j_1,\dots,j_{\ell'})$ (resp. $T^{=}(j_1,\dots,j_{\ell'})$),
to be a vector space spanned by 
all $\Z_{i_1}\cdots \Z_{i_{\ell}}$ with $(i_1,\dots,i_{\ell})>\boldsymbol{j}$ (resp. $(i_1,\dots,i_{\ell})\geq\boldsymbol{j}$).
In order to prove Theorem \ref{mainres}, 
it is enough to show the following seven statements.
\begin{enumerate}
\item[(F1)] For $A,B\in\mathbb{Z}$ with $A>B$, we have $\Z_A\Z_B\in \T(A,B)+R_k$ for any $k\in\mathbb{Z}$.
\item[(F2)] For $A\in\mathbb{Z}$ with $2A\not\in 3\mathbb{Z}$, $\Z_A\Z_A\in \T(A,A)+R_k$ for any $k\in\mathbb{Z}$.
\item[(F3)] For $A\in\mathbb{Z}$ with $2A+1\not\in 3\mathbb{Z}$, $\Z_A\Z_{A+1}\in \T(A,A+1)+R_k$ for any $k\in\mathbb{Z}$.
\item[(F4)] For $i\in\mathbb{Z}$, we have $\Z_{3i}\Z_{3i}\Z_{3i}\in \T(3i,3i,3i)+R_k$ for any $k\in\mathbb{Z}$.
\item[(F5)] For $i\in\mathbb{Z}$, we have $\Z_{3i}\Z_{3i}\Z_{3i+2}\in \T(3i,3i,3i+2)+R_k$ for any $k\in\mathbb{Z}$.
\item[(F6)] For $i\in\mathbb{Z}$, we have $\Z_{3i-2}\Z_{3i}\Z_{3i}\in \T(3i-2,3i,3i)+R_k$ for any $k\in\mathbb{Z}$.
\item[(I)] $V(\Lambda_0+\Lambda_1-\delta)$ and $V(\Lambda_2-2\delta)$ are $\GE$-submodules of $W$.
\end{enumerate}
Here, $R_k$ is the set of expansions of the form (see also ~\cite[(6.19)]{LW3})
\begin{align*}
\xi=\sum_{\ell'\geq 0}\sum_{\boldsymbol{p}=(p_1,\dots,p_{\ell'})\in\mathbb{Z}^{\ell'}}c_{\boldsymbol{p}}\Z_{p_1}\cdots \Z_{p_{\ell'}},
\end{align*}
which may be infinite formal sums,
where $\SUPP_i(\xi)$ is finite for any $i\in\mathbb{Z}$ and $\SUPP_k(\xi)=\emptyset$.
Recall (see ~\cite[(4.8)]{LW3})
that $\SUPP_i(\xi)$ is defined to be the set
\begin{align*}
  \{\boldsymbol{p}=(p_1,\dots,p_{\ell'})\in\mathbb{Z}^{\ell'}\mid c_{\boldsymbol{p}}\ne 0 \textrm{ and }
  p_{\ell''}+\dots+p_{\ell'}\leq i \textrm{ for } 1\leq \ell''\leq\ell'\}.
\end{align*}

The statement (F$a$) (resp. (I)) is proven in \S4.$a+3$ (resp. \S\ref{forseven}) for $1\leq a\leq 6$.

\subsection{Convention}
In the following, for $1\leq a\leq 4$, we denote by $\theta^{(a)}_{A,B}$ the left-hand side minus the right-hand side
of the $a$-th equation in Theorem \ref{mainres}.
Note that $\theta^{(3)}_{A,B}$ and $\theta^{(4)}_{A,B}$ are defined for $A,B\in\mathbb{Z}$ such
that $A+B\not\in 3\mathbb{Z}$.

In the calculation, we use the following explicit values and $c^{(p)}_0=1$ for $1\leq p\leq 5$.
\begin{align*}
\BD &= -24-28\omega+14\omega^3,\quad
\A = -6-8\omega+4\omega^3,\quad
\D = 4-8\omega^2-6\omega^3,\\
\E &= 2-4\omega^2, \quad
\varepsilon(\nu^5(\beta_1),\beta_1) = -52+104\omega^2+90\omega^3,\quad
c^{(6)}_0 = -1-2\omega-\omega^3,\\
c^{(1)}_1 &= \frac{-6-4\omega+2\omega^3}{3},\quad
c^{(2)}_1 = \frac{-4\omega+2\omega^3}{3},\quad
c^{(3)}_1 = \frac{6-4\omega+2\omega^3}{3},\\
c^{(4)}_1 &= \frac{8\omega-4\omega^3}{3},\quad
c^{(5)}_1 = \frac{-6+8\omega-4\omega^3}{3},\quad
c^{(6)}_1 = \frac{4\omega-2\omega^3}{3}.
\end{align*}

\subsection{A proof of (F1)}\label{forone}
For $A>B$, we define
\begin{align*}
\Delta_{A,B}=\ICHI{A,B}\BD(\omega^{4A+9B}+\omega^{9A+4B})-\NI{A,B}\varepsilon(\nu^5(\beta_1),\beta_1)(\omega^{4A+9B}-\omega^{9A+4B}).
\end{align*}
The coefficient $d_{A,B}$ of $\Z_A\Z_B$ in $\Delta_{A,B}$ is given by
\begin{align*}
d_{A,B}=c^{(1)}_0\BD\cdot(\omega^{4A+9B}+\omega^{9A+4B})-c^{(2)}_0\varepsilon(\nu^5(\beta_1),\beta_1)\cdot(\omega^{4A+9B}-\omega^{9A+4B}).
\end{align*}
One can check by case-by-case substitution that
it is non-zero if $A+5\not\equiv B\pmod{12}$.

In the case $A+5\equiv B\pmod{12}$, we consider $-\ICHI{A+1,B-1}+\NI{A+1,B-1}$.
Then, the coefficient of $\Z_A\Z_B$ (resp. $\Z_{A+1}\Z_{B-1}$) is $c^{(2)}_1-c^{(1)}_1=2(\ne0)$ (resp. $c^{(2)}_0-c^{(1)}_0=0$), and that of $\Z'_{A+B}$ is 0 because the following two values are equal.
\begin{align*}
  \varepsilon(\nu^5(\beta_1),\beta_1)(\omega^{4A+9B}-\omega^{9A+4B}) &= \omega^{A+8}(-52+104\omega^2+90\omega^3)(\omega-1),\\
  \BD(\omega^{4A+9B}+\omega^{9A+4B}) &= \omega^{A+8}(-24-28\omega+14\omega^3)(\omega+1).
\end{align*}

\subsection{A proof of (F2)}\label{fortwo}
Consider $\NI{A,A}-\frac{\BD}{\A}\SAN{A,A}$. The coefficient of $\Z_A\Z_A$ is $2(-1-2\omega+\omega^3)$, which is non-zero.

\subsection{A proof of (F3)}\label{forthree}
Consider $c^{(6)}_0(\NI{A,A+1}-\frac{\BD}{\A}\SAN{A,A+1})+(1-\frac{\BD}{\A})\SHI{A,A+1}$. The coefficient of $\Z_A\Z_{A+1}$
is non-zero because of
\begin{align*}
c^{(6)}_0\Big((c^{(2)}_0+c^{(2)}_1)-\frac{\BD}{\A}(c^{(3)}_0+c^{(3)}_1)\Big)+\Big(1-\frac{\BD}{\A}\Big)(c^{(6)}_0-c^{(6)}_1)=8(2+2\omega-\omega^3).
\end{align*}

\subsection{A proof of (F4)}\label{forfour}
Recall \S\ref{forone} and \S\ref{forthree}.
We define
\begin{align*}
  \Delta'_{A,B} &= \Delta_{A,B}/d_{A,B}=\sum_{p\geq 0}c^{(A,B)}_p\Z_{A-p}\Z_{B+p} + d^{(A,B)}\Z_{A+B} + e^{(A,B)},\\
  \Delta'_{C,C+1} &=
\Big(c^{(6)}_0\big(\NI{C,C+1}-\frac{\BD}{\A}\SAN{C,C+1}\big)+\big(1-\frac{\BD}{\A}\big)\SHI{C,C+1}\Big)/(8(2+2\omega-\omega^3))\\
&= \sum_{p\geq 0}c^{(C,C+1)}_p\Z_{C-p}\Z_{C+1+p}+d^{(C,C+1)}\Z_{2C+1} + e^{(C,C+1)}
\end{align*}
for $A, B,C\in\mathbb{Z}$ such that $A>B$ and $A+5\not\equiv B\pmod{12}$, $2C+1\not\in 3\mathbb{Z}$. Note that $c^{(A,B)}_0=c^{(C,C+1)}_0=1$ and
$c^{(C,C+1)}_p$ does not depend on $C$.

Consider $\Delta'_{3i,3i+1}\Z_{3i-1}-\Z_{3i}\Delta'_{3i+1,3i-1}$. It is zero, and has the expansion 
\begin{align*}
  {} &\sum_{1\leq p\leq k-6i}c^{(3i,3i+1)}_p\Z_{3i-p}\Z_{3i+1+p}\Z_{3i-1}\\
  &- (c^{(3i+1,3i-1)}_1\Z_{3i}\Z_{3i}\Z_{3i}
  +\sum_{2\leq r\leq k+1-3i}c^{(3i+1,3i-1)}_r\Z_{3i}\Z_{3i+1-r}\Z_{3i-1+r}).
\end{align*}
modulo $R_k$.
For $p\geq 1$, we have $\Z_{3i+1+p}\Z_{3i-1}\in \TT(3i+p,3i)+R_k$ by \S\ref{forone}.
This, together with $\Z_{3i+1-r}\Z_{3i-1+r}\in \T(3i,3i)$ for $r\geq 2$, implies
\begin{align*}
  \Z_{3i-p}\Z_{3i+1+p}\Z_{3i-1},\Z_{3i}\Z_{3i+1-r}\Z_{3i-1+r}\in \T(3i,3i,3i)+R_k.
\end{align*}
Thus, we have $-c^{(3i+1,3i-1)}_1\Z_{3i}\Z_{3i}\Z_{3i}\in T(3i,3i,3i)+R_k$ and
\begin{align*}
c^{(3i+1,3i-1)}_1=\frac{c^{(1)}_1\BD(1+\omega^{10})-c^{(2)}_1\varepsilon(\nu^5(\beta_1),\beta_1)(1-\omega^{10})}{\BD(1+\omega^{10})-\varepsilon(\nu^5(\beta_1),\beta_1)(1-\omega^{10})}=\frac{\omega(2-\omega^2)}{3}.
\end{align*}

\subsection{A proof of (F5)}\label{forfive}
Recall \S\ref{fortwo} and define
\begin{align*}
  \Delta'_{A,A} = \Big(\NI{A,A}-\frac{\BD}{\A}\SAN{A,A}\Big)/(2(-1-2\omega+\omega^3)) = \sum_{p\geq 0}c^{(A,A)}_p\Z_{A-p}\Z_{A+p}+d^{(A,A)}\Z_{2A} + e^{(A,A)}
\end{align*}
for $A\in\mathbb{Z}$ with $2A\not\in 3\mathbb{Z}$.
Note that $c^{(A,A)}_0=1$ and
$c^{(A,A)}_p$ does not depend on $A$.
By a similar argument as in \S\ref{forfour} applied to $\Delta'_{3i,3i+1}\Z_{3i+1}-\Z_{3i}\Delta'_{3i+1,3i+1}$,
we have $-c_1^{(3i+1,3i+1)}\Z_{3i}\Z_{3i}\Z_{3i+2}\in \T(3i,3i,3i+2)+R_k$ and 
\begin{align*}
c^{(3i+1,3i+1)}_1=\frac{2(c^{(2)}_1-\frac{\BD}{\A}c^{(3)}_1)}{2(-1-2\omega+\omega^3)}=\frac{3+2\omega-\omega^3}{3}.
\end{align*}

\subsection{A proof of (F6)}\label{forsix}
Apply a similar argument as in \S\ref{forfive} to $\Delta'_{3i-1,3i-1}\Z_{3i}-\Z_{3i-1}\Delta'_{3i-1,3i}$, we have
$c^{(3i-1,3i-1)}_1\Z_{3i-2}\Z_{3i}\Z_{3i}\in \T(3i-2,3i,3i)+R_k$.

\subsection{A proof of (I)}\label{forseven}
We assume that readers are familiar with Kashiwara crystals~\cite{Kas2}.
Note that $3\Lambda_0-\alpha_0=\Lambda_0+\Lambda_1-\delta$ and
$3\Lambda_0-2\alpha_0-\alpha_1=\Lambda_2-2\delta$ (see ~\cite[\S6.2, \S12.4]{Kac}).
In order to prove (I), it is enough to show that
$\TOP\otimes\TOP\otimes\KF{0}\TOP$ and
$\TOP\otimes\KF{0}\TOP\otimes\KF{1}\KF{0}\TOP$ are maximal (i.e.,
are sent to $\ZERO$ by $\KE{0},\KE{1},\KE{2}$),
where $\TOP$ is the highest weight element in the highest weight crystal $B(\Lambda_0)$.
This follows from $\KF{1}\KF{0}\TOP\ne\ZERO$, which is easily checked
by explicit realizations such as Kyoto path models~\cite{KKMMNN}, Littelmann path models~\cite{Li2}, etc.

\section{Automatic derivations via the regularly linked sets}\label{automa}


Recall Theorem \ref{mainconj}.
The purpose of this section is to give an automatic derivation of the
generating function $f_{\EL_a}(x,q)$, where we write $f_{\mathcal{C}}(x,q)=\sum_{\lambda\in\mathcal{C}}x^{\ell(\lambda)}q^{|\lambda|}$ for $\mathcal{C}\subseteq\PAR$,
by the regularly linked sets~\cite{TT}, which
generalize Andrews' linked partition ideals~\cite{An0,An1} by finite automata in formal language theory~\cite{Koz}.

\subsection{A survey of ~\cite{TT} for modulus 1}
As in ~\cite[\S3]{TT}, for a non-empty set $\Sigma$, let $\Sigma^{\ast}$ be the set of words $w_1\cdots w_{\ell}$
of finite length of $\Sigma$. By the word concatenation $\cdot$ and the empty word $\EMPTYWORD$, the set $\Sigma^{\ast}$ is regarded as a
free monoid generated by $\Sigma$. For $A,B\subseteq\Sigma^{\ast}$,
we define the sum, the product, and the Kleene star by
\begin{align*}
  A+B=A\cup B,\quad
  AB=\{ww'\mid w\in A,w'\in B \},\quad
  A^{\ast}=\{\EMPTYWORD\}+A+A^2+\cdots.
\end{align*}

\begin{Def}
For a positive integer $\emu$ and subsets $\EFU, \AI\subseteq\PAR\setminus\{\EMPTYPART\}$,
let $\C$ be a subset of $\PAR$ which consists of partitions $\lambda$ such that
\begin{enumerate}
\item $m_j(\lambda)<\emu$ for $j\geq 1$,
\item $\lambda$ does not begin with $\boldsymbol{c}$ for $\boldsymbol{c}\in \AI$, and
\item $\lambda$ does not match $(b_1+k,\dots,b_{p}+k)$ for $(b_1,\dots,b_{p})\in \EFU$ and $k\geq 0$.
\end{enumerate}
\end{Def}

Here, $\EMPTYPART$ is the empty partition (i.e., the partition of 0).
Recall that as in ~\cite[Definition 1.1]{TT}, we say that a partition $\lambda=(\lambda_1,\dots,\lambda_{\ell})$
begins with (resp. matches) $(d_1,\dots,d_r)$ if $\ell\geq r$ and $(\lambda_{\ell-r+1},\cdots,\lambda_{\ell})=(d_1,\dots,d_r)$
(resp. $(\lambda_{i+1},\dots,\lambda_{i+r})=(d_1,\dots,d_r)$  for some $0\leq i\leq\ell-r$).


\begin{Def}
For a non-empty partition $\lambda$ and $\emu\geq 1$, let $\SAT_{\emu}(\lambda)$ (resp. $\SAT'_{\emu}(\lambda)$) be a finite subset of $\PAR$ 
which consists of partitions $\mu$ such that $\mu_1=\lambda_1$, $m_j(\mu)<\emu$ for $j\geq 1$, and $\mu$ matches $\lambda$ (resp. $\mu_1=\lambda_1$, $m_{\lambda_1}(\mu)<\emu$, and $\mu$ begins with $\lambda$).
\end{Def}

\begin{Ex}
  For $\lambda=(5,3,2)$, we have
  $\SAT'_3(\lambda)=\{(5,5,3,2),(5,3,2)\}$. Similarly,
$\SAT_3(\lambda)$ consists of the following partitions.
\begin{align*}
(5,5,3,2,2),
(5,3,2,2),
(5,5,3,2),
(5,3,2),
(5,5,3,2,2,1,1),
(5,3,2,2,1,1),\\
(5,5,3,2,1,1),
(5,3,2,1,1),
(5,5,3,2,2,1),
(5,3,2,2,1),
(5,5,3,2,1),
(5,3,2,1).
\end{align*}
\end{Ex}

Let $\JEI=\{0,\dots,\emu-1\}$ and define a map $\pi:\JEI\to\PAR$ by $\pi(i)=(\underbrace{1,\dots,1}_{i})$ for $i\in\JEI$.
As in ~\cite[Definition 2.5]{TT}, we have an injection
\begin{align*}
\PAI:\SEQ(\JEI,\pi)\hookrightarrow\PAR,
\end{align*}
where $\SEQ(\JEI,\pi)$ is the set of infinite sequences $\boldsymbol{i}=(i_1,i_2,\dots)$ in $\JEI$ (i.e., $i_j\in \JEI$ for $j\geq 1$)
such that $i_j\ne 0$ holds only for finitely many $j$, and
the partition $\mu=\PAI(\boldsymbol{i})$ is characterized by $m_j(\mu)=i_j$ for $j\geq 1$.

By construction, $\PAI$ gives a bijection to $\C$ when restricted to the set
\begin{align*}
\{(i_1,i_2,\dots)\in\SEQ(\JEI,\pi)\mid (i_1,\dots,i_j)\not\in \JEI^{\ast}\cdot\PAII(\SAT_m(\EFU))\cdot\JEI^{\ast}+\PAII(\SAT'_m(\AI))\cdot\JEI^{\ast}\textrm{ for } j\geq 1\},
\end{align*}
where $\PAII(\nu)=m_1(\nu)\cdots m_{\nu_1}(\nu)\in\JEI^{\nu_1}(\subseteq\JEI^{\ast})$
(i.e., $\PAII(\nu)=(m_1(\nu),\dots,m_{\nu_1}(\nu))\in\JEI^{\nu_1}$) for a non-empty partition $\nu$, and $\SAT_m(\EFU)=\bigcup_{\lambda\in\EFU}\SAT_{\emu}(\lambda)$, 
$\SAT'_m(\AI)=\bigcup_{\lambda\in\AI}\SAT'_{\emu}(\lambda)$.

\begin{Rem}
By the notation in ~\cite[Definition 3.8]{TT}, the set above 
is written as $\AVOID(\SEQ(\JEI,\pi),\PAII(\SAT_m(\EFU)),\PAII(\SAT'_m(\AI)))$.
Note also that $\EFU=\emptyset$ (resp. $\AI=\emptyset$) is allowed while we exclude the case
that the empty partition $\EMPTYPART$ belongs to $\EFU$ (resp. $\AI$)
to satisfy the condition ``$X,X'\subseteq\Sigma^{+}$'' in ~\cite[Definition 3.7]{TT}, where $\Sigma^{+}=\Sigma^{\ast}\setminus\{\EMPTYWORD\}$.
\end{Rem}

From now on, we assume $\PAII(\SAT_m(\EFU)), \PAII(\SAT'_m(\AI))\subseteq\JEI^{\ast}$ are
regular~\cite[Definition 3.3]{TT} so that
the right-hand side of \eqref{dfaeq} below is regular~\cite[Proposition 3.5]{TT}.
For example, the assumption is satisfied if $\EFU$ and $\AI$ are finite
(as in \S\ref{finiteap} and \S\ref{finiteap2}).


Let $\EMU=(Q,\JEI,\delta,s,F)$ be a deterministic finite automaton (DFA, for short) such that
(see ~\cite[Definition 3.1]{TT}, ~\cite[Definition 3.2]{TT} and ~\cite[Appendix A]{TT})
\begin{align}
  L(\EMU)=\JEI^{\ast}\cdot\PAII(\SAT_m(\EFU))\cdot\JEI^{\ast}+\PAII(\SAT'_m(\AI))\cdot\JEI^{\ast}.
\label{dfaeq}
\end{align}
Then, ~\cite[Theorem 3.14]{TT} gives a simultaneous $q$-difference equation
\begin{align*}
f_{\C^{(v)}}(x,q)=\sum_{u\in Q\setminus F}\Big(\sum_{\substack{a\in\JEI \\ u=\delta(v,a)}}x^{\ell(\pi(a))}q^{|\pi(a)|}\Big)f_{\C^{(u)}}(xq,q),
\end{align*}
where (see ~\cite[Definition 3.11]{TT} and ~\cite[(3.7)]{TT}) $M_v=(Q,\JEI,\delta,v,F)$ and
\begin{align*}
  \C^{(v)}=\PAI(\{(i_1,i_2,\dots)\in\SEQ(\JEI,\pi)\mid (i_1,\dots,i_j)\not\in L(M_v)\textrm{ for } j\geq 1\}).
\end{align*}
Because of $\C^{(s)}=\C$, it gives a $q$-difference equation for $f_{\C}(x,q)$ in virtue of the Murray-Miller algorithm (see a review in ~\cite[Appendix B]{TT}).


\subsection{An application to $\EL_{3}$}\label{finiteap}
Apply $\emu=3$, $\EFU=F$ (or $\EFU=F\setminus\{(1,1,1)\}$) and $\AI=I_3$ (see Theorem \ref{mainconj}), we get the minimum DFA (unique up to isomorphism)
$\EMU=(Q,\JEI,\delta,\JOO{0},\{\JOO{3}\})$ where 
$Q=\{\JOO{0},\dots,\JOO{10}\}$, $\JEI=\{0,1,2\}$ and
the values $\delta(v,a)$
of the transition function $\delta:Q\times\JEI\to Q$ are displayed as follows.

\begin{center}
\begin{tabular}{r|rrrrrrrrrrr}
$a\backslash v$ & \JO{0} & \JO{1} & \JO{2} & \JO{3} & \JO{4} & \JO{5} & \JO{6} & \JO{7} & \JO{8} & \JO{9} & \JO{10} \\ \hline
\AL{0}          & \JO{1} & \JO{4} & \JO{1} & \JO{3} & \JO{6} & \JO{8} & \JO{6} & \JO{4} & \JO{6} & \JO{5} & \JO{4} \\ 
\AL{1}          & \JO{2} & \JO{5} & \JO{3} & \JO{3} & \JO{7} & \JO{3} & \JO{7} & \JO{2} & \JO{10} & \JO{3} & \JO{3} \\ 
\AL{2}          & \JO{3} & \JO{3} & \JO{3} & \JO{3} & \JO{3} & \JO{3} & \JO{9} & \JO{3} & \JO{3} & \JO{3} & \JO{3}
\end{tabular}
\end{center}

Thus, we have a simultaneous $q$-difference equation
\begin{align*}
  \begin{pmatrix}
    f_{\C^{(\JOO{0})}}(x,q)\\
    f_{\C^{(\JOO{1})}}(x,q)\\
    f_{\C^{(\JOO{2})}}(x,q)\\
    f_{\C^{(\JOO{4})}}(x,q)\\
    f_{\C^{(\JOO{5})}}(x,q)\\
    f_{\C^{(\JOO{6})}}(x,q)\\
    f_{\C^{(\JOO{7})}}(x,q)\\
    f_{\C^{(\JOO{8})}}(x,q)\\
    f_{\C^{(\JOO{9})}}(x,q)\\
    f_{\C^{(\JOO{10})}}(x,q)\\
  \end{pmatrix}
    =
\begin{pmatrix}
0 & 1 & xq & 0 & 0 & 0 & 0 & 0 & 0 & 0 \\ 
0 & 0 & 0 & 1 & xq & 0 & 0 & 0 & 0 & 0 \\ 
0 & 1 & 0 & 0 & 0 & 0 & 0 & 0 & 0 & 0 \\ 
0 & 0 & 0 & 0 & 0 & 1 & xq & 0 & 0 & 0 \\ 
0 & 0 & 0 & 0 & 0 & 0 & 0 & 1 & 0 & 0 \\ 
0 & 0 & 0 & 0 & 0 & 1 & xq & 0 & x^2q^2 & 0 \\ 
0 & 0 & xq & 1 & 0 & 0 & 0 & 0 & 0 & 0 \\ 
0 & 0 & 0 & 0 & 0 & 1 & 0 & 0 & 0 & xq \\ 
0 & 0 & 0 & 0 & 1 & 0 & 0 & 0 & 0 & 0 \\ 
0 & 0 & 0 & 1 & 0 & 0 & 0 & 0 & 0 & 0 
\end{pmatrix} 
\begin{pmatrix}
    f_{\C^{(\JOO{0})}}(xq,q)\\
    f_{\C^{(\JOO{1})}}(xq,q)\\
    f_{\C^{(\JOO{2})}}(xq,q)\\
    f_{\C^{(\JOO{4})}}(xq,q)\\
    f_{\C^{(\JOO{5})}}(xq,q)\\
    f_{\C^{(\JOO{6})}}(xq,q)\\
    f_{\C^{(\JOO{7})}}(xq,q)\\
    f_{\C^{(\JOO{8})}}(xq,q)\\
    f_{\C^{(\JOO{9})}}(xq,q)\\
    f_{\C^{(\JOO{10})}}(xq,q)\\
\end{pmatrix}
\end{align*}


The pseudo-code in ~\cite[Algorithm 1]{TT} stops at the 9-th iteration and
gives a $q$-difference equation for $f_{\C^{(\JOO{0})}}(x,q)$, which is written in Theorem \ref{mainqd} below.

\begin{Thm}\label{mainqd}
For $1\leq a\leq 3$, we have $\sum_{r=0}^{8}p^{(a)}_r(x,q) f_{\EL_a}(xq^r,q)=0$, where
\begin{align*}
p^{(1)}_0 &= (1-xq^{5})(1-x^{2}q^{9}),\quad
p^{(1)}_1=-(1-xq^{5})(1+xq+xq^2-xq^3-xq^4+x^2q^6-x^2q^8-x^2q^9),\\
p^{(1)}_2 &= xq(1-q^{2}-q^{3}+xq^{2}+xq^{3}-2xq^{5}-xq^{6}+xq^{7}+xq^{8}-x^{2}q^{9}+x^{2}q^{11}+x^{2}q^{12}-x^{3}q^{13}),\\
p^{(1)}_3 &= -x^{2}q^{3}(1-q^2)(1+q-q^3+xq^3+xq^4 + xq^5 -xq^6- xq^7 -x^2q^{10}-x^2q^{11}),\\
p^{(1)}_4 &= -x^{2}q^{5}(1+q)(1+xq^4)(1- xq+ xq^4- xq^5- xq^6 +x^2q^8),\\
p^{(1)}_5 &= x^{3}q^{9}(1-q^{2})(1+q- xq- xq^2- xq^3 + xq^4+ xq^5- x^2q^6 - x^2q^7+x^2q^9),\\
p^{(1)}_6 &= -x^{3}q^{10}(1+xq^{4}-xq^{6}-xq^{7}-x^{2}q^{5}-x^{2}q^{6}+2x^{2}q^{8}+x^{2}q^{9}-x^{2}q^{10}-x^{2}q^{11}-x^{3}q^{11}+x^{3}q^{13}+x^{3}q^{14}),\\
p^{(1)}_7 &= x^{4}q^{15}(1-xq^3)(1-q^2-q^3+ xq^3+ xq^4- xq^5 - xq^6+x^2q^{10}),\quad
p^{(1)}_8=x^{4}q^{18}(1-xq^{3})(1-x^{2}q^{7}),\\
p^{(2)}_0 &= (1-xq^{3})(1-xq^{5})(1-xq^{6})(1-x^{2}q^{9}),\\
p^{(2)}_1 &= -(1-xq^{5})(1-xq^{6})(1+xq-xq^3- xq^4- x^2q^4 + x^2q^6 - x^2q^9- x^3q^{10} + x^3q^{11}+x^3q^{12}),\\
p^{(2)}_2 &= -xq^{3}(1-xq^{5})(1-xq^{6})(q-x + x^2q^3- x^2q^5- x^2q^8 +x^3q^9),\\
p^{(2)}_3 &= -x^{2}q^{3}(1-q^2)(1-xq^6)(1+ q^2- xq^5- x^2q^6 - 2x^2q^9- x^2q^{11}+ x^3q^{11}+ x^3q^{13}+x^3q^{14}),\\
p^{(2)}_4 &= -x^{2}q^{6}(1+q)(1+xq^{4})(1-x+xq-xq^2-2xq^4+xq^5-2xq^6+x^2q^3-x^2q^4+x^2q^5+x^2q^6\\
&\quad\quad\quad\quad +x^2q^9+x^2q^{10}+x^2q^{11}-x^3q^8+x^3q^9-x^3q^{10}-2x^3q^{12}+ x^3q^{13}- 2x^3q^{14}+x^4q^{16}),\\
p^{(2)}_5 &= x^{3}q^{8}(1-q^2)(1-xq^2)(1+ q^2+ q^3- xq^3- 2xq^6- xq^8- x^2q^{10}+ x^3q^{13}+x^3q^{15}),\\
p^{(2)}_6 &= -x^{3}q^{11}(1-xq^{2})(1-xq^{3})(1+xq^2- xq^4- xq^7- x^2q^7 +x^3q^{16}),\\
p^{(2)}_7 &= x^{4}q^{17}(1-xq^2)(1-xq^3)(1-q-q^2+xq^2- xq^4+ xq^7-x^2q^7+ x^2q^9+ x^2q^{10}-x^3q^{14}),\\
p^{(2)}_8 &= x^{4}q^{19}(1-xq^{2})(1-xq^{3})(1-xq^{5})(1-x^{2}q^{7}),\\
p^{(3)}_0 &=  (1-x^{2}q^{9})(1-x^{2}q^{10}),\quad
p^{(3)}_1 =  -(1-x^{2}q^{10})(1+xq+xq^2-xq^4-xq^5-x^2q^9),\\
p^{(3)}_2 &=  xq^{2}(1-q^{2}-q^{3}+xq+xq^{2}+xq^{3}-xq^{4}-2xq^{5}-xq^{6}+xq^{8}\\
&\quad\quad\quad\quad +x^{2}q^{8}-x^{2}q^{10}+x^{2}q^{12}+x^{2}q^{13}-x^{3}q^{12}-x^{3}q^{13}+x^{3}q^{15}+x^{3}q^{16}-x^{4}q^{18}),\\
p^{(3)}_3 &= -x^{2}q^{5}(1+q)(1-q^3)(1-q+xq+xq^3-xq^5+x^2q^8-x^2q^9-x^3q^{13}),\\
p^{(3)}_4 &= -x^{2}q^{7}(1+q)(1-xq^2+xq^3+xq^5-xq^6-x^2q^3-x^2q^5+x^2q^6+2x^2q^8\\
&\quad\quad\quad\quad\quad\quad\quad\quad -x^2q^9-x^2q^{10}-x^2q^{11}-x^3q^{10}+x^3q^{11}+x^3q^{13}-x^3q^{14}+x^4q^{16}),\\
p^{(3)}_5 &=  x^{3}q^{11}(1+q)(1-q^3)(1- xq^3+ xq^4- x^2q^4- x^2q^6+ x^2q^8 - x^3q^{11}+x^3q^{12}),\\
p^{(3)}_6 &=  -x^{3}q^{14}(1+xq^{2}+xq^{3}-xq^{5}-xq^{6}-x^{2}q^{6}+x^{2}q^{8}-x^{2}q^{10}-x^{2}q^{11}-x^{3}q^{7}-x^{3}q^{8}\\
&\quad\quad\quad\quad\quad\quad -x^{3}q^{9}+x^{3}q^{10}+2x^{3}q^{11}+x^{3}q^{12}-x^{3}q^{14}-x^{4}q^{14}+x^{4}q^{16}+x^{4}q^{17}),\\
p^{(3)}_7 &= -x^{4}q^{22}(1-x^2q^6)(1-x-xq+xq^3 +xq^4-x^2q^7),\quad
p^{(3)}_8=x^{4}q^{23}(1-x^{2}q^{6})(1-x^{2}q^{7}).
\end{align*}
\end{Thm}

\begin{Rem}
The
$q$-difference equations in Theorem \ref{mainqd} can be obtained by Andrews' linked partition ideals (see ~\cite[Appendix E]{TT}) because $\EFU$ and $\AI$ are finite.
Yet, an approach via the regularly linked sets has an advantage: the minimum forbidden patterns and forbidden prefixes are automatically detected~\cite[Appendix D]{TT}.
In other words, one has a chance to get smaller simultaneous $q$-difference equations.
\end{Rem}

\subsection{A proof of Theorem \ref{mainqd}}\label{finiteap2}
For $\EL_1$ (resp. $\EL_2$), the calculations are similar.
One only needs to perform the automatic calculation (e.g., by computer algebra as in ~\cite[Remark 4.3]{TT}) for the minimum DFA $\EMU$ such that \eqref{dfaeq} for 
$\emu=3$, $\EFU=F$ (or $\EFU=F\setminus\{(1,1,1)\}$) and $\AI=I_1$ (resp. $\AI=I_2$). We omit the details.

\subsection{A proof of Theorem \ref{mainconj}}\label{automaproof}
Note that the case $a=1$ was shown in ~\cite{AvH} by
\begin{align*}
  f_{\EL_1}(x,q) = \sum_{m,k\geq 0}\frac{q^{m^2+3km+4k^2}}{(q;q)_k(q;q)_m}(1-q^{k}+q^{k+m})x^{2k+m}
\end{align*}
(see ~\cite[Theorem 1.7.c)]{AvH}) and by a finite version of a $q$-series identity ~\cite[Proposition 4.3]{AvH}.
We follow a more or less similar line.

We perform an automatic calculation of a $q$-difference equation for each of 
\begin{align*}
  g_2(x,q) &= \sum_{m,k\geq 0}\frac{q^{m^2+3km+4k^2+k+m}}{(q;q)_k(q;q)_m}(1+xq^{3k+1})x^{2k+m},\\
  g_3(x,q) &= \sum_{m,k\geq 0}\frac{q^{m^2+3km+4k^2+2k}}{(q;q)_k(q;q)_m}(1-x^2q^{4m+8k+6})x^{2k+m},
\end{align*}
via a $q$-analog of Sister Celine's technique (see ~\cite[\S7.1]{Tsu}). We see that
it is the same
as that of $f_{\EL_a}(x,q)$ in Theorem \ref{mainqd} for $a=2,3$, which implies
that we have $f_{\EL_a}(x,q)=g_a(x,q)$. 
By taking the limit $n\to\infty$ in ~\cite[Proposition 4.4]{AvH}, we have
\begin{align*}
  f_{\EL_2}(1,q)=\sum_{k\geq 0}\frac{q^{2k^2+k}}{(q;q)_{2k+1}},\quad
  f_{\EL_3}(1,q)=\sum_{k\geq 0}\frac{q^{2k^2+2k}}{(q;q)_{2k+1}}.
\end{align*}
Thanks to Slater's list~\cite[(9), (38)]{Sla}, the right-hand side of the former (resp. latter) equals $1/(q;q^2)_{\infty}$ (resp. $1/(q,q^4,q^6,q^7,q^9,q^{10},q^{12},q^{15};q^{16})_{\infty}$).

\hspace{0mm}

\noindent{\bf Acknowledgments.}
The author thanks Kana Ito and
Matthew Russell for their helpful discussions, especially to Matthew Russell
for bringing the paper~\cite{AvH} to his attention.
This work was supported by the Research Institute for Mathematical
Sciences, an International Joint Usage/Research Center located in Kyoto
University, the TSUBAME3.0 supercomputer at Tokyo Institute of Technology,
JSPS Kakenhi Grant 20K03506, the Inamori Foundation, JST CREST Grant Number JPMJCR2113, Japan and Leading Initiative for Excellent Young Researchers, MEXT, Japan.


\end{document}